\documentclass{article}[12pt]

\usepackage{pstricks}
\usepackage{graphics}
\usepackage{caption}
\usepackage{amsmath}
\usepackage{amsthm}
\usepackage{amsfonts}
\usepackage{amssymb}
\usepackage{amscd}
\usepackage{calc}
\usepackage{graphicx}
\usepackage{xypic}
\usepackage{color}
\usepackage{wasysym}
\usepackage{subcaption}

\newtheorem{thm}{Theorem}[section]

\newtheorem{defi}{Definition}

\newtheorem{quest}{Open Question}

\newcommand{\Nb}{\mathbb{N}}

\newcommand{\UU}{\mathcal{U}}
\newcommand{\KK}{\mathcal{K}}

\newcommand{\BB}{\mathcal{B}}

\newcommand{\GG}{\mathcal{G}}
\newcommand{\HH}{\mathcal{H}}

\newcommand{\tbf}{\textbf}

\newcommand{\bog}[9]{
\xy (0,0)*{}; (0,18)*{} **\dir{-};
(6,0)*{};(6,18)*{} **\dir{-};
(12,0)*{};(12,18)*{} **\dir{-};
(18,0)*{};(18,18)*{} **\dir{-};
(0,0)*{};(18,0)*{} **\dir{-};
(0,6)*{};(18,6)*{} **\dir{-};
(0,12)*{};(18,12)*{} **\dir{-};
(0,18)*{};(18,18)*{} **\dir{-};
(3,3)*{\text{#7}};
(3,9)*{\text{#4}};
(3,15)*{\text{#1}};
(9,3)*{\text{#8}};
(9,9)*{\text{#5}};
(9,15)*{\text{#2}};
(15,3)*{\text{#9}};
(15,9)*{\text{#6}};
(15,15)*{\text{#3}};
\endxy}

\newcommand{\bgraph}{
\xy
(5,5)="13";
(5,15)="12";
(5,25)="11";
(15,5)="23";
(15,15)="22";
(15,25)="21";
(25,5)="33";
(25,15)="32";
(25,25)="31";
"11"*{\bullet}+(-2,2)*{1};"12"*{\bullet}+(-2,0)*{8};"13"*{\bullet}+(-2,-2)*{7};
"21"*{\bullet}+(0,2.5)*{2};"22"*{\bullet}+(1,3)*{9};"23"*{\bullet}+(0,-2.5)*{6};
"31"*{\bullet}+(2,2)*{3};"32"*{\bullet}+(2,0)*{4};"33"*{\bullet}+(2,-2)*{5};
}
\newcommand{\bgB}{
\xy
(3,3)="7";
(3,9)="8";
(3,15)="1";
(9,3)="6";
(9,9)="9";
(9,15)="2";
(15,3)="5";
(15,9)="4";
(15,15)="3";
"1"*{\bullet};"2"*{\bullet};"3"*{\bullet};
"8"*{\bullet};"9"*{\bullet};"4"*{\bullet};
"7"*{\bullet};"6"*{\bullet};"5"*{\bullet};
}

\xyoption{line}

\title{$10$ Questions about Boggle Logic Puzzles}
\author{Jonathan Needleman}
\date{May 21, 2015} 
\begin{document}

\maketitle

\begin{abstract}
Boggle logic puzzles are based on the popular word game Boggle, where you are given  list of words, and  your goal is to recreate a Boggle board.  In this paper we give an overview of known results and then propose a number of problems related to these puzzles.\end{abstract}

\section{Introduction}
 Boggle (\copyright Hasboro inc) is a popular word search game where players compete to find as many words as they can in a $4\times 4$ grid of letters \cite{boggle}.  Boggle logic puzzles on the other hand is the game of Boggle played in reverse.   A list of words are given and you need to recreate the board.  This paper both summarizes the known results and raises $10$ open problems related to these puzzles, many of which are suitable for undergraduate research.  However, before we explore these questions, we first give an overview the rules of Boggle.

For a word to be valid it must be at least three letters long, and  consecutive letters in the word must be adjacent on the board, either horizontally, vertically or  diagonally.  Furthermore, each box on the board cannot be used more than once in forming a word.  For instance the word MOMATH is on the board because of the path in figure \ref{BOARD}.
\begin{figure}[h!]
\[
\xy (0,0)*{}; (0,24)*{} **\dir{-};
(6,0)*{};(6,24)*{} **\dir{-};
(12,0)*{};(12,24)*{} **\dir{-};
(18,0)*{};(18,24)*{} **\dir{-};
(24,0)*{};(24,24)*{} **\dir{-};
(0,0)*{};(24,0)*{} **\dir{-};
(0,6)*{};(24,6)*{} **\dir{-};
(0,12)*{};(24,12)*{} **\dir{-};
(0,18)*{};(24,18)*{} **\dir{-};
(0,24)*{};(24,24)*{} **\dir{-};
(3,3)="41";
(3,9)="31";
(3,15)="21";
(3,21)="11";
(9,3)="42";
(9,9)="32";
(9,15)="22";
(9,21)="12";
(15,3)="43";
(15,9)="33";
(15,15)="23";
(15,21)="13";
(21,3)="44";
(21,9)="34";
(21,15)="24";
(21,21)="14";
"11"*{G};"12"*{O};"13"*{M};"14"*{C};
"21"*{M};"22"*{T};"23"*{I};"24"*{B};
"31"*{H};"32"*{A};"33"*{N};"34"*{S};
"41"*{U};"42"*{E};"43"*{A};"44"*{W};
"13"*{}; "12"*{} **\dir{-};"12"*{}; "21"*{} **\dir{-};"21"*{}; "32"*{} **\dir{-};
"32"*{}; "22"*{} **\dir{-};"22"*{}; "31"*{} **\dir{-};
\endxy\]
\caption{A Boggle board}\label{BOARD}
\end{figure}
However, BIBS is not on the board, because there is only one ``B'' on the board, and HUG is not on the board because the ``U'' and ``G'' are not adjacent.

A Boggle logic puzzle is a list of words that can be found in a unique board (up to rotation and reflection).  To get a feel for this  try the following puzzle:

$$\textbf{act, ape, ate, cop, end, old}$$

The goal is to fill in a $3\times 3$ board so that each of the six words appear according to the rules of Boggle.  In $2007$ Mark Zegarelli published a book \cite{Zeg} of Boggle logic puzzles.  All of his puzzles are on a full $4\times 4$ board, but they all have at least one letter filled in.  This is most likely to guarantee uniqueness (there are no symmetries in his puzzles).

This article is an attempt to summarize what is known about the mathematics of Boggle logic puzzles, and then give the reader a number of open questions that touch upon a wide variety of fields.  However, to do this, we need to translate the problem into mathematical language.  We will also make some assumptions to simplify the initial problems.

\section{The mathematical language of Boggle}
\subsection{Basic notation}
In order to make the problems more tractable, we assume every letter on a given board is unique.  This will prevent us from having to look at many different cases depending on which letters repeat and where they are placed.  In section \ref{sec:repeat} we will address some questions that pertain to puzzles with repeated letters.

For an $n\times n$ board, instead of using letters, we will use the numbers $1 \ldots n^2$.  By doing so we ignore English spelling of words, and let any sequence of adjacent, non-repeating numbers be a word.  We denote words with hyphens between ``letters'' such as $2-4-5$.  The set of all $n\times n$ boards of this form are denoted by $\BB_n$.

To this point most research has focused on $3\times 3$ boards, so we fix a standard $3 \times 3$ board $\bf{S}$ to be as in figure \ref{fig:main}
\begin{figure}[!htbp]
  \begin{subfigure}[b]{.6\linewidth}
    \centering
    $\bog{1}{2}{3}{8}{9}{4}{7}{6}{5}$
    \caption{}
    \label{fig:main}
  \end{subfigure}%
  \begin{subfigure}[b]{.2\linewidth}
    \centering
   $\bog{2}{1}{8}{9}{3}{4}{7}{6}{5}$
    \caption{}
    \label{fig:perm}
  \end{subfigure}%
  \caption{The standard board $\bf{S}$, and a permutation.}
  \label{fig:boards}
\end{figure}

It is necessary to talk about other boards besides $\bf{S}$.  To do this we use $\bf{S}$ and permutations.  Let $S_9$ be permutations of the numbers $[1 \ldots 9]$.  These permutations to refer to other boards in the following way:  if $g\in S_9$ then we let $[g]$ be the board $g\bf{S}$.  For example, the element $(12)(359)\in S_9$ gives the   board $[(12)(359)]$ seen in figure \ref{fig:perm} by switching the $1$ and the $2$, and also replacing the $3$ with a $5$, the $5$ with a $9$ and a $9$ with a $3$.

In order to reference words, let $W(\bf{B})$ be the set of all words  in the board $\bf{B}$ that are at least two letters long.  Since there is no way to distinguish two different boards that contain exactly the same  words (such as a reflection of the board), we say two boards $\bf{B}, \bf{B}'\in \BB_n$ are equivalent  if and only if $W(\tbf{B})=W(\tbf{B}')$, and write this as $\tbf{B}\cong \tbf{B}'$.  It turns out $\tbf{B}\cong \tbf{B}'$ if and only if the boards are a rotation or a reflection of each other.

This terminology allows us to mathematically define Boggle logic puzzles.

\begin{defi}
An $n\times n$ {\bf Boggle logic puzzle} is a list of words $P$ so that there is $\tbf{B}\in\BB_n$ with $P\subset W(\tbf{B})$ and whenever there is another $\tbf{B'}\in\BB_n$ with $P\subset W(\tbf{B}')$ then $\tbf{B}\cong \tbf{B}'$.
\end{defi}
All this is saying is that the list of words (or puzzle) $P$ can be found in one and only one board up to equivalence. 

\subsection{Boggle and graphs}

While it is good to have a formal definition of these puzzles, we still need a language in which to prove various properties of these puzzles.  Graph theory provides a great a language to do this.  This is because Boggle logic puzzles are entirely about which letters are adjacent to each other, and graph theory is the mathematics of adjacency.

We can view the Boggle board as a graph by thinking of boxes on the board as vertices, with edges placed between adjacent boxes.  For a $n\times n$ board we call this graph $\KK_n$ because it is the well known $n\times n$ king's graph.  Furthermore, we can view $\KK_n$ as a labeled graph by labeling the vertices with the letters of the standard board.  For example, the labeled $\KK_3$ is shown in figure \ref{fig:K3}.  In the future we will suppress the labels to prevent clutter.
\begin{figure}[!htbp]
        \centering
                    \[\bgraph
                        "11"*{};"12"*{}**\dir{-};"12"*{};"13"*{}**\dir{-};"11"*{};"21"*{}**\dir{-};"21"*{};"31"*{}**\dir{-};
                        "31"*{};"32"*{}**\dir{-};"32"*{};"33"*{}**\dir{-};"13"*{};"23"*{}**\dir{-};"23"*{};"33"*{}**\dir{-};
                        "21"*{};"22"*{}**\dir{-};"22"*{};"23"*{}**\dir{-};"12"*{};"22"*{}**\dir{-};"22"*{};"32"*{}**\dir{-};
                        "11"*{};"22"*{}**\dir{-};"12"*{};"21"*{}**\dir{-};"22"*{};"31"*{}**\dir{-};"21"*{};"32"*{}**\dir{-};
                        "13"*{};"22"*{}**\dir{-};"22"*{};"33"*{}**\dir{-};"32"*{};"23"*{}**\dir{-};"12"*{};"23"*{}**\dir{-};
                    \endxy\]
                    \caption{Labeled $\KK_3$}
                    \label{fig:K3}
            \end{figure}

\begin{figure}[!htbp]
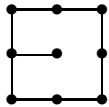
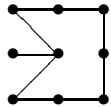

                \begin{subfigure}[b]{.5\linewidth}
                    \centering
                    \[\bgB
                        "1";"2"**\dir{-};"3";"2"**\dir{-};"4";"3"**\dir{-};"4";"5"**\dir{-};
                        "5";"6"**\dir{-};"6";"7"**\dir{-};"7";"8"**\dir{-};"8";"9"**\dir{-};"1";"8"**\dir{-};
                    \endxy\]
                    \caption{$\GG(P)$}
                    \label{fig:GP}
               \end{subfigure}%
               \begin{subfigure}[b]{.4\linewidth}
                    \centering
                    \[\bgB
                        "1";"2"**\dir{-};"3";"2"**\dir{-};"4";"3"**\dir{-};"4";"5"**\dir{-};
                       "5";"6"**\dir{-};"6";"7"**\dir{-};"9";"8"**\dir{-};"1";"9"**\dir{-};"7";"9"**\dir{-};
                    \endxy\]
                    \caption{A non-labeling subgraph}
                    \label{fig:nonlab}
               \end{subfigure}%
               \caption{Non-labeling subgraphs}
            \end{figure}

In addition to viewing the board as a graph, we can also interpret the puzzles themselves as graphs.  For a standard list of words $P\subset W(\bf{S})$ we view the ``letters'' in the words of $P$ as vertices, and place edges between them if the letters are adjacent in some word of $P$.  This graph $\GG(P)$ is called the {\bf adjacency graph} of $P$.  For instance if $P=\{9-1-2, 2-3-4-5, 5-6-7-9, 7-9-8\}$ then $\GG(P)$ would be as in figure \ref{fig:GP}.  The vertices are placed in the same relative position as in the labeled $K_3$ graph for consistency.

The language of graph theory allows us to determine when two puzzles $P_1, P_2\subset W(\textbf{B})$ are really the same puzzle.   If $\GG(P_1)$ and $\GG(P_2)$ are isomorphic as graphs then we say $P_1$ is \textbf{equivalent} to $P_2$ and write $P_1\cong P_2$.  In non-technical language this is  saying both puzzles $P_1$ and $P_2$ provide the exact same information about which letters are adjacent to others.

Solving a Boggle logic puzzle $P$ becomes an exercise of figuring out how to view an adjacency graph $\GG(P)$ as a subgraph of $\KK_n$.  If you can do  this, the labeling from $\GG(P)$ can then be used to solve the puzzle.    With this in mind, we introduce what we believe to be a novel definition in graph theory.

\begin{defi}
Let $\HH$ be a subgraph of $\GG$.  Then $\HH$ is a \textbf{labeling subgraph} of $\GG$ if whenever any subgraph $\HH'$  of $\GG$ is isomorphic to $\HH$  by $\phi$, then the isomorphism $\phi$ induces an automorphism on $\GG$.
\end{defi}
This definition is rather technical, but the idea behind it is that any labeling of $\mathcal{H}$ induces a unique, up to symmetry, labeling of $\GG$.    To better understand what this definitions means, let us look at a non-example for $\KK_3$.  The graph in figure \ref{fig:GP} is not a labeling subgraph of $\KK_3$.  This is because of the graph in figure \ref{fig:nonlab}.  These two graphs are isomorphic subgraphs of $\KK_3$ by exchanging the left hand side vertex (labeled $8$ in $\KK_3$) with the center vertex (labeled $9$ in $\KK_3$).  However, this exchange in not an isomorphism of $\KK_3$ since $8$ is adjacent to $5$ vertices while $9$ is adjacent to $8$ vertices in $K_3$.

This definition is designed so that $P$ is a puzzle if and only if $\GG(P)$ is a labeling subgraph of $\KK_n$.  In 1972 Entringer and Erd\~os introduced a similar idea called a unique subgraph \cite{unique}.  The idea is that a subgraph $\HH$ of $\GG$ is unique if and only if there is no other subgraph isomorphic to $\HH$.  For a puzzle $P$  to have a  unique solution it may seem $\GG(P)$ only needs to be a unique subgraph of $\KK_n$ up to automorphism of $\KK_n$.  However, this is not strong enough as seen by the graph $\UU$ in figure \ref{fig:uniq}.  All subgraphs $\KK_3$ that are isomorphic to $\UU$ come from automorphism of $\KK_3$, that is, they are just rotations and reflections of $\UU$.  However, this is not a labeling subgraph because $\UU$ is isomorphic to itself by a map the swaps the two upper corners ($1$ and $3$).  This map is not an automorphism of $\KK_3$.  In terms of puzzles this means that both the boards in figures \ref{fig:sol1} and \ref{fig:sol2} would be solutions to any puzzle $P$ with $\GG(P)=\UU$.

\begin{figure}[!htbp]
    \begin{subfigure}[b]{.3\linewidth}
        \centering
                    \[\bgB
                        "1";"2"**\dir{-};"3";"2"**\dir{-};"4";"9"**\dir{-};"8";"9"**\dir{-};"5";"6"**\dir{-};"6";"7"**\dir{-};
                        "1";"9"**\dir{-};"2";"9"**\dir{-};"3";"9"**\dir{-};"5";"9"**\dir{-};"6";"9"**\dir{-};"7";"9"**\dir{-};
                        "2";"4"**\dir{-};"2";"8"**\dir{-};"4";"6"**\dir{-};"6";"8"**\dir{-};"8";"7"**\dir{-};"4";"5"**\dir{-};
                    \endxy\]
                    \caption{$\UU$}
                    \label{fig:uniq}
    \end{subfigure}
    \begin{subfigure}[b]{.3\linewidth}
        \centering
                    $\bog{1}{2}{3}{8}{9}{4}{7}{6}{5}$
                    \caption{A solution to $\UU$}
                    \label{fig:sol1}
    \end{subfigure}\begin{subfigure}[b]{.4\linewidth}
        \centering
                    $\bog{3}{2}{1}{8}{9}{4}{7}{6}{5}$
                    \caption{Another solution}
                    \label{fig:sol2}
    \end{subfigure}
    \caption{Unique but non-labeling}
\end{figure}

\section{Known results: The extremes}\label{sec_max}
There are two natural questions concerning the extremal behavior of Boggle logic puzzles.  The first is, ``What are the fewest words needed to create a puzzle''?  The second is ``How long does a list of words need to be in order to guarantee that you can uniquely recreate the board''?

For a $3\times 3$ board the first question is answered in \cite{bogmin} with the following theorem:

\begin{thm}\label{thm:min}
Any $3\times 3$ Boggle logic puzzle with no repeated letters and using only three-letter words must contain at least $6$ different words.
\end{thm}

This result is readily extended to longer word lengths.

The second question has not been answered as satisfactory, but has the following partial answer is presented in \cite{bogmax}.

\begin{thm}\label{thm:max}
For a $3\times 3$ Boggle board with no letters repeated, one needs $137$ different (out of $160$ possible) three-letter words to guarantee a unique solution.  For four-letter words, $377$ words are needed out of $496$ possible words.
\end{thm}

\subsection{Minimal Solutions}

The paper \cite{bogmin} shows any $3\times 3$ Boggle logic puzzle must be made of at least $6$ three letter words.  This paper studies two-letter words instead of looking at three-letter words.  In fact, the proof of  theorem \ref{thm:min} only proves that at least $11$ two-letter words are needed to create a puzzle.   This is strong enough to claim that at least $6$ three-letter words are needed.   However, every  puzzle found so far uses at least $12$ two-letter words.  This leads to our first open question.

\begin{quest}\label{Q:min11}
Are there any $3\times 3$ puzzles with nine distinct letters  that have only $11$ two-letter words?
\end{quest}
Our conjecture is that the answer to this question is no.  Some modifications to the proof of theorem \ref{thm:min} is probably all that is necessary, making this problem accessible to undergraduates.

The next question, however, is more involved and introduces the idea of minimal puzzles.  A \textbf{minimal puzzle} $P$  is a puzzle where if you break up a word into two smaller words (possibly a one letter word) then $P$ is no longer a puzzle.  Equivalently, removing an edge from $\GG(P)$ yields a  graph that is not a labeling subgraph.  Essentially, a minimal puzzle is a puzzle where all the information given in the puzzle is needed to solve the puzzle.  This begs the following question:

\begin{quest}\label{Q:mincount}
How many inequivalent minimal $n\times n$ puzzles are there?
\end{quest}
From a puzzle maker's point of view this is asking about how robust Boggle logic puzzles are as a source of puzzles.  Will they be making the same puzzles over and over just with different words?  We expect that $3\times 3$ boards are not a great source of puzzles, but $4\times 4$ and larger boards probably will yield a large number of puzzles.

Once question \ref{Q:min11} has been answered there is another version of question \ref{Q:mincount} to answer for $3\times 3$  boards.  How many minimal puzzles are there with either $11$ or $12$ two-letter words (depending on the answer to question \ref{Q:min11})?  This question has the potential to be easier because all puzzles with either $11$ or $12$ words would be minimal, so one only needs to check how many words are in the puzzle.  For a $3\times 3$ boards there are either ${20}\choose{11}$ or ${20}\choose{12}$ word lists that may be minimal, and so these questions can probably be answered by computer.  However a $4\times 4$ board has $42$ edges instead of $20$ so the number of subgraphs of a fixed number of edges is likely larger than a computer can handle, so it would be ideal to have analytic solutions.

\subsection{Maximal Solutions}\label{sec:max}
Another natural question about Boggle logic puzzles comes from the game of Boggle itself.  Imagine playing a round of Boggle and at the end of the round you have a long list of words.  Is there an easy way to determine if you can recreate the board from your list?  The first approach you may take is to see if there is a number $N$ where, if you have that many words, you are guaranteed to be able to recreate the board.  Theorem \ref{thm:max} from \cite{bogmax} is a such a result for $3\times 3$ boards with three- and four-letter words.

The idea of the proof is relatively simple:  find a board $\textbf{B}$ that maximizes the number of words in common with the standard board $\textbf{S}$.  One more than the number of common words is the cut-off value $N$ to guarantee a unique solution.  For both three- and four-letter words this board is $\textbf{[}(1 \ 3)\textbf{]}$ as seen in figure \ref{fig:sol2}.  However, the proof for three-letter words and four-letter are essentially done as independent proofs.  This leads to some open questions which should be accessible to undergraduates.
\begin{quest}\label{Q:uproof}
Does the  board $\textbf{[}(1 \ 3)\textbf{]}$ have the most $k$-letter words in common with $\textbf{S}$ of any $3\times 3$ board, and if so is there a uniform proof?
\end{quest}
The $3\times 3$ board appears to have some degeneracy that does not appear in larger boards.  In particular, the board $\textbf{[}(1 \ 3)\textbf{]}$  , which swaps adjacent corners, is special in the $3\times 3$ case due to the fact that both $1$ and $3$ are both adjacent to a common side letter $2$.  Since this is not the case for $n \times n$ boards with $n\geq 4$ swapping adjacent corners probably will not maximize common words.  Instead we provide the following conjecture in the form of an open question.

\begin{quest}\label{Q:max}
Given an $n\times n$ board $\textbf{B}$ with $n\geq 4$, does the board which comes from swapping a corner letter with an adjacent side letter maximize the number of $k$-letter words in common with $\textbf{B}$?
\end{quest}

\section{Open Questions}
In this section we introduce a number of open questions unrelated to any known result.

\subsection{Probabilistic Questions}
In section \ref{sec:max} we saw that it takes a large number of words to guarantee a unique solution, probably more than one could expect to find in any actual game of Boggle.  Instead one can ask probabilistic questions in hopes of obtaining answers that may be more in line with what one may see in a game of Boggle.  With this in mind we ask two questions, one about the mean, and one about the median.

\begin{quest}\label{Q:mean}
For a given $n\times n$ board what is the expected number of $k$-letter words needed to recreate the board?
\end{quest}

\begin{quest}\label{Q:med}
For a given $n\times n$ board what is the number of $k$-letter words needed so that there is a $50\%$ chance the board can be recreated?
\end{quest}
The naive brute force methods to answer these questions quickly run into problems.  Let us just think about these problems in the simplest case, a $3\times 3$ board $\textbf{B}$ and two-letter words.

For question \ref{Q:mean}, one approach is to take a random word $w_1\in W(\textbf{B})$, and determine the average number of boards that contain $w_1$.  Then choose another word $w_2$ at random and determine how many boards are expected to include both $w_1$ and $w_2$.  Continue adding words until you expect only one board.  For one word this computation is doable, but once you add more words the dependence structure becomes quite complicated.  The number of boards containing words $w_1,\ldots w_m$ depends greatly on the adjacency structure of the words chosen.  This makes a direct computation of the mean number of boards containing these words difficult.

On the other hand, an answer to question \ref{Q:med} is likely related to the answer to question \ref{Q:mincount}.  If one can determine how many length $m$ lists of two-letter words have a unique solution then one essentially has an answer to question \ref{Q:med} since the total number of lists is relatively easy to count.

\subsection{Computational Complexity}
It is natural to ask about the computational complexity of solving a Boggle logic puzzle.  Up to this point we have been thinking of puzzles on a $n\times n$ boggle board.  For a given $n\times n$ puzzle $P$ this amounts to determining how the adjacency graph $\GG(P)$ can be viewed as a subgraph of $\KK_n$.   Instead of playing on a standard Boggle board the puzzle can actually be done on any graph.  From this perspective solving Boggle logic puzzles is related to solving the subgraph isomorphism problem, a classic NP-complete problem.

The subgraph isomorphism problem asks if a graph $\HH$ is isomorphic to a subgraph of $\GG$.  This problem can be shown to be NP-complete by showing it is equivalent to the clique problem (is a complete graph a subgraph) or by the Hamiltonian cycle problem (does a graph have a Hamiltonian cycle), both of which are known to be NP-complete \cite{NPcomplete}.

Returning to Boggle puzzles on standard boards we observe that the  graph $\KK_n$ does not contain a complete graph of order $5$, and it contains a Hamiltonian cycle for all $n$.   This makes it unclear if solving Boggle logic puzzles is NP-complete, and so we raise the following question.

\begin{quest}\label{Q:solveB}
What is the computational complexity of solving $n\times n$ Boggle logic puzzles?
\end{quest}

A related and potentially more difficult question is:

\begin{quest}\label{Q:createB}
What is the computational complexity of creating $n\times n$ Boggle logic puzzles?
\end{quest}

Question \ref{Q:createB} is asking ``Is a graph $\HH$ a labeling subgraph of $\KK_n$''?  For the subgraph isomorphism problem one wants the existence of an isomorphism from $\HH$ to a subgraph of $\KK_n$.  For the labeling subgraph problem one wants to know there exists {\it exactly} eight isomorphisms from $\HH$ to subgraphs of $\KK_n$.  Since this is a counting problem creating Boggle logic puzzles is at least as difficult as solving them.

There is no reason to limit oneself to the computational complexity of solving Boggle logic puzzles instead of asking  about labeling subgraphs of general graphs.  The only difference is that if $\HH$ is a labeling subgraph of $\GG$ one needs to show the number of isomorphisms from $\HH$ to subgraphs of $\GG$ is equal to the number of automorphisms of $\GG$.

\subsection{Repeated Letters}\label{sec:repeat}
To this point we have only discussed problems that assume all the letters on the board are different.  Of course, all those questions can be asked again assuming that some letters are allowed to repeat, but as we will see this does not necessarily make sense in all contexts.  Instead, we introduce some questions about  which boards with repeated letters  are worthy of study for Boggle logic puzzles.

Let us look at all $3\times 3$ boards with eight $1$s and one $2$ for letters.  It is not hard to check that every board with these letters must contain the same words.   In other words, it does not matter where the $2$ is placed.  When all letters are different the only equivalence between boards occur from the symmetries of a square, but in this case all boards are equivalent.  This means it does not make sense to study Boggle logic puzzles with these letters since there is no list of words that will distinguish boards.

Before moving on we introduce a way to discuss which letters are being used in a given board.  Let $\lambda=(\lambda_1,\ldots \lambda_m)$ be a partition of $n^2$.  That is $\sum \lambda_i = n^2$ with $\lambda_i\in \Nb$ and $\lambda_i\geq \lambda_{i+1}$.  We will say a $n\times n$ board is type $\lambda$ if it contains $\lambda_i$ copies of the letter $i$.  For instance, boards of type $\lambda=(8,1)$ are all the $3\times 3$ boards with eight $1$s, and a single $2$.  Whereas  boards of the form  $\lambda=(1,1,1,1,1,1,1,1,1)$ (which we will abbreviate as $\lambda(1^9)$) are all boards with nine distinct letters.  We then let $\BB(\lambda)$ be all boards of type $\lambda$.

We have seen boards of type $(8,1)$ are not interesting to study as Boggle logic puzzles, so the question becomes which boards actually are interesting to study.  We will say board $\textbf{B}$ is {\bf solvable} if all their equivalent boards come from symmetries of a square.  We know if all the letters are different this is the case, but this isn't always true when there are repeated letters.
\begin{quest}\label{Q:intboard}
Let  $\lambda$ be a partition of $n^2$.  Which $\textbf{B}\in \BB(\lambda)$ are solvable?
\end{quest}
For $\lambda=(1^{n^2})$ all boards are solvable, but this is not an interesting result because no letters are repeated.  Are there any other solvable boards?  Yes!  In fact all boards of type $(2, 1^7)$ (boards with a single letter repeated once) are solvable.  There is no published proof of this fact so we briefly outline the ideas here.

First, boards with with two adjacent $1$s are inequivalent to boards with non-adjacent $1$s.  For boards with two adjacent $1$s look at how many words can be of the form $1-x-1$ and $1-1-x$.  This pair of numbers is distinct for each placement of two adjacent $1$s.  For non-adjacent $1$s, checking how many words are of the form $1-x-1$ is not quite enough as there are two different boards that have two words of that form, but then those boards can be seen to be inequivalent.

This leads to one final open question.

\begin{quest}\label{Q:allboard}
For which  $\lambda$, a partition of $n^2$, are all $\textbf{B}\in \BB(\lambda)$ are solvable?
\end{quest}

\section{Conclusions}
Boggle logic puzzles are fun puzzles that have rich mathematics at their core.  The fact that the computational complexity of these puzzles is related to the subgraph isomorphism problem shows how the puzzle is related to an area of active research.  On the other hand, the puzzles give rise to other problems in graph theory that are accessible to people with a wide range of backgrounds.

The study  of Boggle logic puzzles is very young, and the problems listed in this article is only intended to be a starting place for further research.  We hope that in the next edition of this journal many of these questions will be solved and many more new problems proposed.


\begin{thebibliography}{9}
\bibitem{unique}
     Entringer, R. C. and Erd\~os, P, {\it  On the number of unique subgraphs of a graph}. J. Combinatorial Theory B 13, 1972, 112-115.
\bibitem{NPcomplete}
     Garey, Michael R.  and Johnson, David S.  {\it Computers and Intractability: A Guide to the Theory of NP-Completeness}. W.H. Freeman.  1979.
\bibitem{boggle}
    Hasbro Incorporated.  Boggle Game, {\it Rule Book}, 1997.
\bibitem{bogmin}
    Needleman, Jonathan. {\it Boggle Logic Puzzles:  Minimal Solution}.  College Math Journal,
Vol. 44, No. 4 September 2013, 293-299.
\bibitem{bogmax}
    Needleman, Jonathan and Porrino, Lauren. {\it Boggle Logic Puzzles: Maximal non-solutions}. Accepted, Pi Mu Epsilon.
\bibitem{Zeg}
    Zegarelli, Mark. {\it Sit and Solve BOGGLE Logic Puzzles} Puzzlewright 2007.

\end{thebibliography}
\end{document}